\documentclass[12pt]{article}

\usepackage{latexsym}
\usepackage{amsfonts}
\usepackage{amssymb}
\usepackage{amsmath}
\usepackage{mathrsfs}
\usepackage{hyperref}

\newcommand{\be}{\begin{equation}}
\newcommand{\ee}{\end{equation}}
\newcommand{\bea}{\begin{eqnarray}}
\newcommand{\eea}{\end{eqnarray}}
\newcommand{\bean}{\begin{eqnarray*}}
\newcommand{\eean}{\end{eqnarray*}}
\newcommand{\brray}{\begin{array}}
\newcommand{\erray}{\end{array}}

\newtheorem{dfn}{Definition}[section]
\newtheorem{thm}[dfn]{Theorem}
\newtheorem{lmma}[dfn]{Lemma}
\newtheorem{ppsn}[dfn]{Proposition}
\newtheorem{crlre}[dfn]{Corollary}
\newtheorem{xmpl}[dfn]{Example}
\newtheorem{rmrk}[dfn]{Remark}

\newcommand{\bdfn}{\begin{dfn}\rm}
\newcommand{\bthm}{\begin{thm}}
\newcommand{\blmma}{\begin{lmma}}
\newcommand{\bppsn}{\begin{ppsn}}
\newcommand{\bcrlre}{\begin{crlre}}
\newcommand{\bxmpl}{\begin{xmpl}}
\newcommand{\brmrk}{\begin{rmrk}\rm}

\newcommand{\edfn}{\end{dfn}}
\newcommand{\ethm}{\end{thm}}
\newcommand{\elmma}{\end{lmma}}
\newcommand{\eppsn}{\end{ppsn}}
\newcommand{\ecrlre}{\end{crlre}}
\newcommand{\exmpl}{\end{xmpl}}
\newcommand{\ermrk}{\end{rmrk}}

\newcommand{\clh}{\mathcal{H}}

\newcommand{\clk}{\mathcal{K}}

\author{S.P. Murugan and S. Sundar}
\title{$E_0^{P}$-semigroups and product systems}

\begin{document}
\maketitle
\begin{abstract}
Let $P$ be a closed convex cone in $\mathbb{R}^{n}$. Assume that $P$ is spanning i.e. $P-P=\mathbb{R}^{n}$ and pointed i.e. $P \cap -P=\{0\}$. Let  $\alpha:=\{\alpha_{x}:x \in P\}$ be a $\sigma$-weakly continuous family of 
unital normal endomorphisms on $B(\clh)$. Denote the "product system" associated to $\alpha$ by $\mathcal{E}_{\alpha}$. We show that $\mathcal{E}_{\alpha}$ is a concrete product system and $\alpha$, up to cocycle 
conjugacy, can be recovered completely from $\mathcal{E}_{\alpha}$.

\end{abstract}

\noindent {\bf AMS Classification No. :} {Primary 46L55; Secondary 46L99.}  \\
{\textbf{Keywords.}} $E_0^{P}$-semigroups, Product systems, Laplace transform, Dual Cone.

\section{Introduction}
The theory of $E_{0}$-semigroups, invented by Powers and further developed extensively by Arveson, is now more than 30 years old.  We refer the reader to the monograph \cite{Arveson}  for a comprehensive view of the subject.   In short, in  Arveson's programme of $E_0$-semigroups, one  attempts to study and possibly classify actions of the semigroup of non-negative reals, $\mathbb{R}_{+}$,
on the algebra of bounded operators on a separable Hilbert space or more generally on von-Neumann algebras.  Arveson proposes a classification of $E_0$-semigroups into three broad types namely 
Type I, II and III of which only Type I $E_0$-semigroups are completely classified. We should mention here that the classification problem in the Type II and Type III cases are far from over. 

This paper, along with the paper \cite{Anbu}, attempts to enlarge Arveson's programme to more general semigroups i.e. attempts to study unital actions of more general semigroups on $B(\clh)$, 
where $\clh$ is an infinite dimensional separable Hilbert space.  The motivation for enlarging Arveson's programme stems from the fact that most of the  concrete examples in the theory of $E_0$-semigroups 
arise from isometric representations of $\mathbb{R}_{+}$ or perturbations of them.  But isometric representations of more general semigroups and in particular the $C^{*}$-algebras associated to them
have been  studied extensively for more than three decades now. We refer the reader to \cite{Li13} and the references therein for discrete semigroup $C^{*}$-algebras and to \cite{SundarJFA} and the references therein for topological semigroup $C^{*}$-algebras. 
We believe that extending Arveson's programme to more general semigroups requires borrowing of ideas from semigroup $C^{*}$-algebras. It is this interaction that interests us and it is the raison d'$\hat{e}$tre of this paper.
This paper is technical in nature and discussion of concrete examples are taken up in  \cite{Anbu}.

Now let us briefly discuss the contents of this paper.  One of the main tools in the theory of $E_0$-semigroups is the notion of a product system. Arveson associates to each $E_0$-semigroup a measurable field 
of Hilbert spaces having an associative product structure.  The associated field of Hilbert spaces is called the  product system associated to the given $E_0$-semigroup. It turns out that the associated product
system is in fact a complete invariant. We extend this result in the context of $E_{0}^{P}$-semigroups (Defn. \ref{E-semigroups}) where $P$ is a closed convex cone of $\mathbb{R}^{n}$ which is spanning and pointed.
This is achieved in Thm.\ref{Main theorem}.  The major technical hurdle that needs overcoming is to prove that  there is no distinction between measurable $E_{0}^{P}$-semigroups and
continuous $E_{0}^{P}$-semigroups. This is proved in  Prop.\ref{technical lemma}. The proof of Prop.\ref{technical lemma} utilises ideas from convex analysis, the notion of a dual cone and the Laplace transform of a 
function over a cone.  We should mention here that the proof of Prop.\ref{technical lemma} is heavily inspired by the ideas proposed in \cite{Faraut}, \cite{Renault_Muhly}, \cite{Nica_WienerHopf} and in \cite{Hilgert_Neeb}.

\section{Preliminaries}
Here we collect the basic definitions that are  needed to read this paper. 
Throughout this paper, the letter $P$ stands for a closed convex cone in $\mathbb{R}^{n}$. We  assume that $P$ is spanning i.e. $P-P=\mathbb{R}^{n}$ and  pointed i.e. $P \cap -P=\{0\}$. 
We denote the interior of $P$ by $\Omega$.  Note that $\Omega$ is an ideal in $P$ in the sense that $P + \Omega \subset \Omega$. For $x,y \in \mathbb{R}^{n}$, we write $x \leq y$ and $x<y$ if $y-x \in P$ and $y-x \in \Omega$ respectively.  We need the fact that $\Omega$ is dense in $P$. This is standard in convex analysis. Nevertheless we include the proof but postpone it till  next section. 

 The letter $\clh$ stands 
for an infinite dimensional separable Hilbert space. We denote the von-Neumann algebra of bounded operators on $\clh$ by $B(\clh)$ and we denote its predual, the trace class operators by
$\mathcal{L}^{1}(\clh)$. For a von-Neumann algebra $M$, its predual will be denoted by $M_{*}$. For $\rho \in M_{*}$ and $A \in M$, we write $\rho(A)$ as $\langle \rho | A \rangle$. We consider only von-Neumann algebras with separable predual. %We also treat $M$ as a measurable space where the measurable structure is given by the $\sigma$-algebra generated by $\sigma$-weakly closed subsets of $M$. 

\begin{dfn}
\label{E-semigroups}
Let $M$ be a von-Neumann algebra with  separable predual $M_{*}$. By  a  $\sigma$-weakly continuous/measurable $E^{P}$-semigroup on $M$, we mean a family $\alpha=\{\alpha_{x}\}_{x \in P}$ of normal endomorphisms on $M$, indexed by $P$, such that the following conditions hold: 
\begin{enumerate}
\item[(1)] For $x,y \in P$, $\alpha_{x}\circ \alpha_{y}=\alpha_{x+y}$.
\item[(2)] For $\rho \in M_{*}$ and $A \in M$, the map $P \ni x \to \rho(\alpha_{x}(A)) \in \mathbb{C}$ is continuous/measurable.
\item[(3)]  $\alpha_{0}$ is the identity map.
\end{enumerate}
If in addition, if  $\alpha_{x}(1)=1$ for every $x \in P$, we call the family $\alpha=\{\alpha_{x}\}_{x \in P}$ a $\sigma$-weakly continuous/measurable $E_{0}^{P}$-semigroup on $M$. 
\end{dfn}
\begin{rmrk}
 We will see later  in Prop.\ref{technical lemma} that there is no distinction between $\sigma$-weakly continuous $E^{P}$-semigroups and $\sigma$-weakly measurable $E^{P}$-semigroups. Till then, by an $E^{P}$-semigroup on $M$, we mean a $\sigma$-weakly continuous one. But we retain the adjective $\sigma$-weakly measurable while referring to a $\sigma$-weakly measurable $E^{P}$-semigroup. 
 \end{rmrk}
 
 \begin{rmrk}
If $M=B(\clh)$, we simply refer an $E_{0}^{P}$-semigroup on $M$ as an $E_{0}^{P}$-semigroup.
If $P=\mathbb{R}_{+}$, we recover Arveson's notion of $E_0$-semigroups. In Arveson's programme of $E_0$-semigroups, one  attempts to study and possibly classify  $E_0$-semigroups up to cocycle conjugacy. This paper
along with the paper \cite{Anbu} are the first steps towards enlarging Arveson's programme to more general semigroups. In this paper, we restrict ourselves to closed convex cones. However the basic definitions go through
when $P$ is a closed subsemigroup of a locally compact group containing the identity element.  
\end{rmrk}

\begin{rmrk}
Let $\alpha:=\{\alpha_{x}\}_{x \in \Omega}$ be a semigroup of  normal $*$-endomorphisms of $M$, indexed by $\Omega$. We say that $\alpha$ is a $\sigma$-weakly continuous/measurable $E^{\Omega}$-semigroup on $M$ if for every $\rho \in M_{*}$ and $A \in M$, the map $\Omega \ni x \to \rho(\alpha_{x}(A)) \in \mathbb{C}$ is $\sigma$-weakly continuous/measurable.  If $\alpha:=\{\alpha_{x}\}_{x \in \Omega}$ is an $E^{\Omega}$-semigroup with $\alpha_{x}$ being unital for each $x \in \Omega$, then we call $\alpha$ an $E_{0}^{\Omega}$-semigroup on $M$. We will see in Prop.\ref{technical lemma} that every  $E_{0}^{\Omega}$-semigroup admits a unique   $E_{0}^{P}$-semigroup extension.
\end{rmrk}

\textbf{Cocycle conjugacy:}  Let $\alpha:=\{\alpha_{x}:x \in P\}$ be an $E^{P}$-semigroup on a von-Neumann algebra $M$. A $\sigma$-weakly continuous family of unitaries in $M$, $\{u_{x}:x \in P\}$, is called an $\alpha$-cocycle if $u_{x}\alpha_{x}(u_y)=u_{x+y}$.  Let $u:=\{u_{x}:x \in P\}$ be an $\alpha$-cocycle. For $x \in P$, set $\beta_{x}=Ad(u_x) \circ \alpha_{x}$. Then $\beta:=\{\beta_{x}:x \in P\}$ is an $E^{P}$-semigroup on $M$ and is called a cocycle perturbation of $\alpha$. Note that $\beta$ is an $E_{0}^{P}$-semigroup on $M$ if and only if $\alpha$ is an $E_{0}^{P}$-semigroup on $M$. Let $\alpha:=\{\alpha_{x}: x \in P\}$ and $\beta:=\{\beta_{x}:x \in P\}$ be  $E^{P}$-semigroups on $M$ and $N$ respectively. We say that $\alpha$ and $\beta$ are cocycle conjugate if there exists a normal isomorphism $\theta:M \to N$ such that $\{\theta \circ \alpha_{x} \circ \theta^{-1}: x \in P\}$ is a cocycle perturbation of $\beta$.

One of the important tools in Arveson's  classification  programme of $E_0$-semigroups is the notion of a product system. Arveson associates to each $E_0$-semigroup a "product system" which turns out to be a complete invariant [See \cite{Arveson}]. 
The goal of this paper is to enlarge the  aforementioned invariant to $E_{0}^{P}$-semigroups. %\newline

Let $\alpha:=\{\alpha_{x}\}_{x \in P}$ be an $E_{0}^{P}$-semigroup on $B(\clh)$. For $x \in P$, let \[
E_{x}:=\{T \in B(\clh): \alpha_{x}(A)T=TA ~\forall A \in B(\clh)\} .\]
The vector subspace $E_{x}$ is called the  intertwining space of $\alpha_{x}$. Let $T,S \in E_{x}$ be given. Note that $T^{*}S$ commutes with every bounded operator on $\clh$. Hence $T^{*}S$ is a scalar. We denote this scalar by $<T,S>$. With respect to the inner product $<,>$, $E_{x}$ is a Hilbert space. Also the Hilbert space norm on $E_{x}$ coincides with the operator norm. Moreover observe that if $T \in E_{x}$, $S \in E_{y}$ then $TS \in E_{x+y}$.  Note that the linear span of $\{TS: T \in E_{x}, S \in E_{y}\}$ is dense in $E_{x+y}$. The proof of this fact follows exactly as in the $1$-dimensional case i.e. when $P=[0,\infty)$.  Thus we omit the proof.

Observe that for $x \in P$, the intertwining space $E_{x}$ is either $1$-dimensional or of infinite dimension. For, $\{\alpha_{tx}: t \geq 0\}$ is an $E_0$-semigroup and an application of  Theorem 2.4.7 in \cite{Arveson} yields the desired conclusion.

\begin{lmma}
Let $\alpha:=\{\alpha_{x}: x \in P\}$ be an $E_0^{P}$-semigroup on $B(\clh)$. The following are equivalent. 
\begin{enumerate}
\item[(1)] There exists $x \in \Omega$ such that $E_{x}$ is $1$-dimensional.
\item[(2)] For every $x \in P$,  $E_{x}$ is $1$-dimensional.
\item[(3)] There exists $x \in \Omega$ such that $\alpha_{x}$ is an automorphism. 
\item[(4)] For every $x \in P$, $\alpha_{x}$ is an automorphism.
\end{enumerate}
\end{lmma}
\textit{Proof:} The structure theorem for a single endomorphism (Prop. 2.1.1, \cite{Arveson}) implies that $(1)$ is equivalent to $(3)$  and $(2)$ is equivalent to $(4)$. Thus, it is enough to prove that $(3)$ implies $(4)$. Let $s \in \Omega$ be such that $\alpha_{s}$ is an automorphism. For $x \in P$, let $M_{x}:=\alpha_{x}(B(\clh))$. Note that if $x \leq y$ then $M_{y} \subset M_{x}$. Consider a point $x \in P$. Since $s \in \Omega$ and $\Omega$ is open, it follows that there exists $n \in \mathbb{N}$ such that $y:=s-\frac{x}{n} \in \Omega$. Thus $\frac{x}{n} \leq s$. As a consequence, we have $B(\clh)=M_{s} \subset M_{\frac{x}{n}}$. This implies that $\alpha_{\frac{x}{n}}$ is onto and hence an automorphism.  Now It follows immediately that  $\alpha_{x}$ is an automorphism. This completes the proof.

\begin{rmrk}
$E_{0}^{P}$-semigroups for which each endomorphism is an automorphism are completely classified by $\mathbb{T}$-valued $2$-cocycles on $\mathbb{R}^{n}$ (\cite{Anbu}). In view of this, we do not consider such degenerate cases. 
\end{rmrk}
In this paper, by an $E_{0}^{P}$-semigroup, we mean a family $\alpha:=\{\alpha_{x}\}_{x \in P}$ of endomorphisms on $B(\clh)$ as in Defn. \ref{E-semigroups} such that for every $x \in \Omega$, $\alpha_{x}$ is \textbf{not onto}.

The following definition of a concrete product system is inspired by Defn.2.4.2 of \cite{Arveson}. We simply replace the open interval $(0,\infty)$ by $\Omega$, the interior of $P$. We recall the definition in full for the reader's sake. We treat $B(\clh)$ as a measurable space where the measurable structure on $B(\clh)$ is given by the $\sigma$-algebra generated by $\sigma$-weakly closed subsets of $B(\clh)$. Let $p:\Omega \times B(\clh) \to \Omega$ be the projection onto the first co-ordinate. 

\begin{dfn}
\label{concrete product systems}
Let $\mathcal{E} \subset \Omega \times B(\clh)$ be a measurable subset. For $x \in \Omega$, let \[\mathcal{E}(x):=\{T \in B(\clh): (x,T) \in \mathcal{E}\}.\] 
We say that $\mathcal{E}$ is a concrete product system over $\Omega$ if the following conditions are satisfied:
\begin{enumerate}
\item[(1)] The map $p:\mathcal{E} \to \Omega$ is onto.
\item[(2)] For $x \in \Omega$, $\mathcal{E}(x)$ is a vector subspace of $B(\clh)$. Also for $x \in \Omega$, $T,S \in \mathcal{E}(x)$, $T^{*}S$ is a scalar which we denote by $<T,S>$. With respect to the inner product $<,>$, $\mathcal{E}(x)$ is a Hilbert space. 
\item[(3)] For $x,y \in \Omega$, the linear span of $\{TS: T \in \mathcal{E}(x), S \in \mathcal{E}(y)\}$ is dense in $\mathcal{E}(x+y)$. 
\item[(4)] There exists a sequence $\{V_1, V_2, \cdots \}$ of measurable maps from $\Omega$ to $B(\clh)$ such that for every $x \in \Omega$, $\{V_{1}(x),V_{2}(x),\cdots,\}$ is  an orthonormal basis for $\mathcal{E}(x)$. 
\end{enumerate}
\end{dfn}

\begin{rmrk}
Two remarks are in order.
\begin{itemize}
%\item Condition $(4)$ of Defn. \ref{concrete product systems} can also be rephrased by saying that the measurable field of Hilbert spaces $(\mathcal{E}(x))_{x \in \Omega}$ is isomorphic to the trivial field $\Omega \times \clh_{0}$ for an infinite dimensional separable Hilbert space $\clh_{0}$. 
\item Condition $(4)$ of Defn. \ref{concrete product systems} is equivalent to the fact that for every $x \in \Omega$, $\mathcal{E}(x)$ is infinite dimensional and there exists a sequence $\{V_1,V_2,\cdots\}$ of measurable maps from $\Omega \to B(\clh)$ such that for every $x \in \Omega$, $\{V_1(x),V_2(x),\cdots,\}$ is total in $\mathcal{E}(x)$. For a proof of this equivalence, we refer the reader to Prop.7.27 of \cite{Folland}.

\item Let $\mathcal{E} \subset \Omega \times B(\clh)$ and $\mathcal{F} \subset \Omega \times B(\clk)$ be concrete product systems.  We say that $\mathcal{E}$ is isomorphic to $\mathcal{F}$ if for every $x \in \Omega$, there exists an unitary operator $U_{x}:\mathcal{E}(x) \to \mathcal{F}(x)$ such that the maps \[\mathcal{E} \ni (x,T) \to (x,U_{x}(T)) \in \mathcal{F}~, ~\mathcal{F} \ni (x,T) \to (x,U_{x}^{*}(T)) \in \mathcal{E}\] are measurable and \[
U_{x+y}(TS)=U_{x}(T)U_y(S) \]
for $x,y \in \Omega$, $T \in \mathcal{E}(x)$ and $S \in \mathcal{E}(y)$. 

\end{itemize}
\end{rmrk}

 Let $\alpha:=\{\alpha_{x}\}_{x \in P}$ be an $E_{0}^{P}$-semigroup on $B(\clh)$. Define 
\[
\mathcal{E}_{\alpha}:=\{(x,T) \in \Omega \times B(\clh): \alpha_{x}(A)T=TA ~~\forall A \in B(\clh)\}.\]

The goal of this paper is to prove the following  theorem.

\begin{thm} 
\label{Main theorem}
With the foregoing notations, we have the following.
\begin{enumerate}
\item[(1)] Let $\alpha:=\{\alpha_{x}\}_{x \in P}$ be an $E_{0}^{P}$-semigroup on $B(\clh)$. Then $\mathcal{E}_{\alpha}$ is a concrete product system. We call $\mathcal{E}_{\alpha}$ the product system associated to $\alpha$. 
\item[(2)] Suppose $\alpha$ and $\beta$ are  two $E_{0}^{P}$-semigroups on $B(\clh)$. Then $\alpha$ and $\beta$ are cocycle conjugate if and only if $\mathcal{E}_{\alpha}$ and $\mathcal{E}_{\beta}$ are isomorphic as concrete product systems.
\end{enumerate}
\end{thm}

\section{Convex analysis}
In this section, we collect the necessary technical tools from convex analysis that we need to prove Theorem \ref{Main theorem}. 
We refer the reader to Chapter 1 of \cite{Faraut} for proofs. Let us first recall the notion of a dual cone. The dual cone of $P$, denoted $P^{*}$, is defined
as 
\[
P^{*}:=\{u \in \mathbb{R}^{n}: \langle u|x \rangle \geq 0 ,~\forall x \in P\}.
\]
Then $P^{*}$ is also a closed, convex, spanning and a pointed cone. It is well known that, the second dual of $P$, i.e. $(P^{*})^*$ is $P$ itself. We denote the interior of $P^{*}$ by $\Omega^{*}$.  The interior of $P^{*}$ is given by 
\[
\Omega^{*}:=\{a \in \mathbb{R}^{n}: \langle a | x \rangle >0, ~\forall x \in P\backslash\{0\}\}.
\]

The facts that we need about closed convex cones are summarised below.

\begin{lmma}
\label{dominated convergence theorem}
We have the following. 
\begin{enumerate}
\item[(1)] The interior of $P$, $\Omega$ is dense in $P$.
\item[(2)] The boundary of $P$, i.e. $P \backslash \Omega$ has Lebesgue measure zero.
\item[(3)] Let $(a_{k})$ be a sequence in $\mathbb{R}^{n}$ such that $(a_{k})\to a$.  Then $1_{P+a_{k}} \to 1_{P+a}$ a.e. Also $1_{\Omega+a_k} \to 1_{\Omega+a}$ a.e. 
\end{enumerate}
\end{lmma}
\textit{Proof.} We leave it to the reader to convince himself that to prove $(1)$, it is enough to show that $\Omega$ is non-empty. The fact that $\Omega$ is non-empty is a direct consequence of Prop. I.1.4 of \cite{Faraut}.  For a proof of $(2)$, we refer the reader to Lemma 4.1.\ of \cite{Jean_Sundar}. Note that if $a_k\to a$ then $1_{P+a_{k}}$ converges pointwise to $1_{P+a}$ on the complement of the boundary of $P+a$.  The almost convergence of $1_{\Omega+a_{k}}$ to $1_{\Omega+a}$ follows from the fact that the boundary of $P$ has Lebesgue measure zero. This completes the proof. \hfill $\Box$

We now come to the main technical ingredient i.e. the Laplace transform over a cone. 

\begin{dfn}
Let $f:\Omega \to \mathbb{C}$ be bounded and measurable. The Laplace transform of $f$, denoted $L(f)$, is a function on $\Omega^{*}$ defined by the following formula: For $a \in \Omega^{*}$, \[L(f)(a):=\int_{\Omega} e^{-\langle a|x\rangle} f(x)dx.\]
\end{dfn}
\begin{rmrk}
The Laplace tranform of the constant function $'1'$ is usually called the characteristic function of the cone $\Omega$ \cite{Faraut}. 
\end{rmrk}
Note that the Laplace transform is well-defined. To see this, use Lemma I.1.5 of \cite{Faraut} to observe that for $a \in \Omega^{*}$,  there exists $M>0$ and $k$ large such that for $x \in P$, \[e^{-<a|x>} \leq  \frac{M}{(1+||x||^{2})^{k}}.\] Another application of Lemma I.1.5 of \cite{Faraut} implies that the Laplace transform is continuous. 

\begin{ppsn}
\label{Laplace transform}
Let $f:\Omega \to \mathbb{C}$ be bounded and measurable. If $L(f)(a)=0$ for every $a \in \Omega^{*}$ then $f=0$ a.e.
\end{ppsn}
\textit{Proof:} First assume that $f \in L^{1}(\Omega)$. We consider $f$ as a function on $L^{1}(P)$ by declaring its value outside $\Omega$ to be zero. A simple application of Lemma I.1.5. of \cite{Faraut} and Stone-Weierstrass theorem implies that the linear span of $\{e^{-\langle a|x \rangle}: a \in \Omega^{*}\}$ is dense in $C_{0}(P)$. As a consequence, it follows that $\int_{P} \phi(x)f(x) dx=0$ for every $\phi \in C_{0}(P)$. But $f \in L^{1}(P)$. Consequently,  $f=0$ a.e. 

Now fix $a_{0} \in \Omega$ be given. Define $g:\Omega \to \mathbb{C}$ by $g(x):=e^{-\langle a_{0}|x \rangle}f(x)$. Clearly $g \in L^{1}(\Omega)$. The given hypothesis implies that $L(g)=0$. By what we have proved so far, it follows that $g=0$ a.e. which is equivalent to saying that $f=0$ a.e. This completes the proof. \hfill $\Box$

\begin{rmrk}
We should remark that $P$ is pointed is needed to show that the exponentials $\{e^{-\langle a|x\rangle}: a \in \Omega^{*}\}$ separate points of $P$.  

\end{rmrk} 

\section{Measurability versus continuity}
This section forms the soul of this paper. Let $M$ be a von-Neumann algebra with separable predual $M_{*}$. 
In what follows, we treat $M$ and $M_{*}$ as  measurable spaces where the measurable structures are given by the $\sigma$-algebra generated by the weak$^*$-closed subsets of $M$ and weakly closed subsets of $M_{*}$ respectively. 
We also consider $M$ and $M_{*}$ as topological spaces where the topologies considered are the $\sigma$-weak topology on $M$ and the weak topology on $M_*$.

Let $\alpha$ be a normal endomorphism on $M$. By duality, there exists a linear map $\beta:M_{*} \to M_{*}$ such that \[\langle \beta(T),A \rangle = \langle T, \alpha(A) \rangle\] for $T \in M_{*}$ and $A \in M$. We say that $\alpha$ and $\beta$ are dual to each other. 

%Let $\alpha:=\{\alpha_{x}\}_{x \in \Omega}$ be a semigroup of  normal endomorphisms on $M$.  
%By duality, there exists a semigroup of linear maps $\{\beta_{x}\}_{x \in \Omega}$  on $M_{*}$ such that for $T \in M_{*}$, $A \in M$, and $x \in \Omega$, 
%\begin{equation}
%\ langle \beta_{x}(T), A \rangle = \langle T, \alpha_{x}(A)\rangle. 
%\end{equation}
 The main technical hurdle that needs overcoming is Proposition \ref{technical lemma}. When $P=\mathbb{R}_{+}$, this is exactly Proposition 2.3.1 of \cite{Arveson}.  %We simply replace the open interval $(0,\infty)$ by $\Omega$. 
  We start with a little lemma.
 
 \begin{lmma}
 \label{boundary continuity}
 Let $\alpha:=\{\alpha_{x}:x \in P\}$ be a semigroup of normal, unital, $*$-endomorphisms on $M$. Assume that $\alpha$ is faithful in the sense that $\bigcap_{x \in \Omega}\ker (\alpha_{x})=\{0\}$. Suppose that $\alpha$ restricted to $\Omega$ i.e. $\{\alpha_{x}: x \in \Omega\}$ is an $E_{0}^{\Omega}$-semigroup on $M$. Then $\alpha$ is an $E_{0}^{P}$-semigroup on $M$.
  \end{lmma}
\textit{Proof.} The faithfulness of $\alpha$ implies that $\alpha_{0}$ is the identity map. Now the only  thing that demands verification is the $\sigma$-weak continuity of $\alpha$. Let $A \in M$ and $(x_k)$ be a sequence in $P$ converging to $x \in P$. We claim that $\alpha_{x_k}(A)$ converges $\sigma$-weakly to $\alpha_{x}(A)$.  The $\sigma$-weak compactness of $\{T \in M: ||T|| \leq ||A||\}$ allows us to assume that $\alpha_{x_k}(A) \to B$ for some $B \in M$. We claim that $B=\alpha_{x}(A)$.

 Let $s \in \Omega$ be given. 
 Note that $x_k+s$ is a sequence in $\Omega$ that converges to $x+s$. The given hypothesis implies  that $\alpha_{x_k+s}(A) \to \alpha_{s+x}(A)=\alpha_{s}(\alpha_{x}(A))$.  On the other hand, $\alpha_{x_k+s}(A)=\alpha_{s}(\alpha_{x_k}(A)) \to \alpha_{s}(B)$. Thus $\alpha_{s}(B)=\alpha_{s}(\alpha_{x}(A))$ for every $s \in \Omega$.  The fact that $\alpha$ is faithful implies that  implies that $B=\alpha_{x}(A)$. This completes the proof. \hfill $\Box$

\begin{ppsn}
\label{technical lemma}
Let $\alpha:=\{\alpha_{x}: x \in \Omega\}$ be a $\sigma$-weakly measurable $E^{\Omega}$-semigroup on $M$.  Assume that $\alpha$ is faithful i.e. $ \bigcap_{x \in \Omega} Ker(\alpha_{x})=\{0\}$. 
%Suppose that $\alpha$, restricted to $\Omega$, is $\sigma$-weakly measurable i.e. the map $\Omega \ni x \to \langle T  , \alpha_{x}(A) \rangle \in \mathbb{C}$ is measurable  for $A \in M$ and $T \in M_{*}$.
 
 \begin{enumerate}
 \item[(i)] Then $\alpha$ is an $E^{\Omega}$-semigroup on $M$.
 \item[(ii)] Suppose $\alpha$ is an $E_{0}^{\Omega}$-semigroup. Then  there exists a unique $E_{0}^{P}$-semigroup extension $\{\widetilde{\alpha}_{x}: x \in P\}$ on $M$ such that $\widetilde{\alpha}_{x}=\alpha_{x}$ for every $x \in \Omega$. 
 \end{enumerate}
   %$\sigma$-weakly continuous i.e.  the map $P \ni x \to \langle T, \alpha_{x}(A)\rangle \in \mathbb{C}$ is continuous for $A \in M$ and $T \in M_{*}$.
\end{ppsn}
\textit{Proof:} Let $\{\beta_{x}: x \in \Omega \}$ be the dual semigroup of linear maps on $M_{*}$. Note that $\{\beta_{x}: x \in \Omega\}$ is norm bounded. For $a \in \Omega^{*}$ and $T \in M_{*}$, let \[
T(a):=\int_{\Omega} e^{-\langle a | x\rangle} \beta_{x}(T) dx.
\]
We claim the following. \begin{enumerate}
\item[(a)] For $T \in M_{*}$ and $a \in \Omega^{*}$, the map $\Omega \ni x \to \beta_{x}(T(a)) \in M_{*}$ is weakly continuous. 
\item[(b)] The linear span of $\{T(a): T \in M_{*}, a \in \Omega^{*}\}$ is norm dense in $M_{*}$.
\item[(c)] Let $(x_k)$ be a sequence in $\Omega$ such that $x_{k} \to x $ where $x \in P$. Then for every $T \in M_{*}$, $\beta_{x_k}(T)$ converges weakly.  Moreover $\displaystyle \lim_{n \to \infty}\beta_{x_k}(T)$ is independent of the chosen sequence $(x_k)$. 
\end{enumerate}
Since $\{\beta_{x}:x \in P\}$ is norm bounded, we see that once $(a)$ and $(b)$ are established, it is immediate that $(i)$ holds. Let $T \in M_{*}$ and $a \in \Omega^{*}$ be given. Consider a sequence  $(x_{k})$  in $\Omega$ such that $x_{k} \to x \in \Omega$. Calculate as follows to observe that
\begin{align*}
\beta_{x_k}(T(a))&=\int_{\Omega}e^{-\langle a | y \rangle}\beta_{x_k+y}(T)dy \\
                            &= \int_{\mathbb{R}^{n}}e^{-\langle a |u-x_k\rangle}1_{\Omega+x_k}(u)\beta_{u}(T)du 
  \end{align*}
A simple application of dominated convergence theorem, together with Lemma \ref{dominated convergence theorem} and Lemma I.1.5 of \cite{Faraut}, implies that $\beta_{x_k}(T(a))$ converges weakly to $\beta_{x}(T(a))$.  

 Suppose $(b)$ is not true. Then there exists a non-zero $A \in M$ such that for $T \in M_{*}$ and $a \in \Omega^{*}$, \[\langle T(a), A \rangle=\int_{\Omega} e^{-\langle a | x \rangle}\langle \beta_{x}(T), A \rangle dx =0.\] Now Prop. \ref{Laplace transform}
and the fact that $M_{*}$ is separable implies that there exists a subset $E \subset \Omega$ of measure zero such that if $x \notin E$ then $\langle \beta_{x}(T), A \rangle = \langle T, \alpha_{x}(A)\rangle =0$ for every $T \in M_{*}$. 
Thus for $x \in \Omega \backslash E$, $\alpha_{x}(A)=0$. Let $y \in \Omega$ be given. Observe that $y-\Omega$ is an open subset of $\mathbb{R}^{n}$ containing $0$. Since $\Omega$ is dense in $P$, it follows that $(y-\Omega) \cap \Omega$ is a 
non-empty open subset of $\Omega$. Since $E$ has measure zero, there exists $x \in \Omega \cap (y-\Omega) \cap E^{c}$. Write $y=x+z$ with $z \in \Omega$. Since $\alpha_{x}(A)=0$, it follows that $\alpha_{y}(A)=0$. But $\alpha$ is faithful. Hence $A=0$ which is a contradiction. This completes the proof of $(b)$. Now $(i)$ is immediate. 

To prove $(c)$, since $\{\beta_{x}\}_{x \in \Omega}$ is norm bounded, it is enough to show that for $T \in M_{*}$ and $a \in \Omega^{*}$, $\beta_{x_k}(T(a))$ converges weakly. This follows immediately by the integral expression of $T(a)$ and $(3)$ of 
Lemma \ref{dominated convergence theorem}. In fact, the limit of $\beta_{x_k}(T(a))$ as $n$ approches $\infty$ is \[\int_{\mathbb{R}^{n}}e^{-\langle a |u-x\rangle}1_{\Omega+x}(u)\beta_{u}(T)du.\] Now $(c)$ is immediate. For $x \in P$ and $T \in M_{*}$, define $\widetilde{\beta}_{x}(T):=\lim_{n \to \infty}\beta_{x_k}(T)$ where $(x_k)$ is any sequence in $\Omega$ converging to $x$. Here the  limit is taken in the weak sense. Clearly $\widetilde{\beta}_{x}=\beta_{x}$ for $x \in \Omega$.

Let us turn our attention to $(ii)$. Fix $x \in P$ and $A \in M$.

\textit{Claim:} Suppose $(x_{k})$, $(y_k)$ are sequences in $\Omega$  and $B,C \in M$ are such that $x_{k} \to x$, $y_k \to x$, $\alpha_{x_k}(A) \to B$ and $\alpha_{y_k}(A) \to C$. Then $B=C$. 
Let $s \in \Omega$ be given. Note that $x_{k}+s \to x+s$ and $y_k+s \to x+s$. But $x+s \in \Omega$. Observe that 
\[
\alpha_{s+x_k}(A)=\alpha_{s}(\alpha_{x_{k}}(A)).
 \]                            
Taking limits, we obtain $\alpha_{s+x}(A)=\alpha_{s}(B)$. Similarly we obtain $\alpha_{s+x}(A)=\alpha_{s}(C)$. Thus $\alpha_{s}(B-C)=0$ for every $s \in \Omega$. The faithfulness of $\alpha$ implies that $B=C$. 
This proves our claim. 

Now consider a sequence $(x_k)$  in $\Omega$ such that $(x_k) \to x$. We claim that $\alpha_{x_k}(A)$ converges $\sigma$-weakly.
 Note that $\{T \in M: ||T|| \leq ||A||\}$ is $\sigma$-weak compact. Hence it suffices to show that any convergent subsequence of $(\alpha_{x_k}(A))$ converges to the same limit which follows from what is proved in the preceeding paragraph. Moreover the argument in the preceeding paragraph implies that the limit of $(\alpha_{x_k}(A))$ is independent of the chosen sequence $(x_k)$.  Let $\widetilde{\alpha}_{x}(A):=\lim_{n \to \infty} \alpha_{x_k}(A)$. Clearly $\widetilde{\alpha}_{x}=\alpha_{x}$ if $x \in \Omega$. 
 
 Let $x \in P$ be given. We claim that $\widetilde{\alpha}_{x}(U)$ is a unitary if $U$ is a unitary. Let $U \in M$ be a unitary and set $V:=\widetilde{\alpha}_{x}(U)$. Choose a sequence $(x_k)$ in $\Omega$ such that $x_k \to x$ and let $s \in \Omega$ be given. Now note that $\alpha_{x_k+s}(U) \to \alpha_{x+s}(U)$ which is a unitary. But $\alpha_{x_k+s}(U)=\alpha_{s}(\alpha_{x_k}(U))$ converges, by definition, to $\alpha_{s}(\widetilde{\alpha}_{x}(U))=\alpha_{s}(V)$. Thus $\alpha_{s}(V^{*}V)=\alpha_{s}(VV^{*})=1$. Since $\alpha$ is faithful, it follows that $V^{*}V=VV^{*}=1$. Consequently it follows that  $\widetilde{\alpha}_{x}(U)$ is a unitary.

 We leave it to the reader to verify that $\widetilde{\alpha}_{x}$ is a unital $*$-endomorphism of $M$.  We only indicate that this uses the fact that the linear span of unitaries of $M$ is $M$.
  To see, that $\widetilde{\alpha}_{x}$ is normal, observe that for $T \in M_{*}$ and $A \in M$, $\langle \widetilde{\beta}_{x}(T) | A \rangle = \langle T | \widetilde{\alpha}_{x}(A) \rangle$. Now it is immediate that $\widetilde{\alpha}_{x}$ is normal. To see that $\{\widetilde{\alpha}_{x}\}_{x \in P}$ is a semigroup, fix $x,y \in P$ and choose sequences $(x_{k})$ and $(y_k)$ in $\Omega$ such that $x_{k} \to x$ and $y_k \to y$. Fix $A \in M$. Fix $m$ and let $n$ tends to infinity in the equation $\alpha_{x_k+y_{m}}(A)=\alpha_{x_k} \circ \alpha_{y_{m}}(A)$ to obtain the equation $\alpha_{x+y_{m}}(A)=\widetilde{\alpha}_{x}\circ \alpha_{y_m}(A)$. Now let $m \to \infty$ to obtain the equation $\widetilde{\alpha}_{x+y}(A)=\widetilde{\alpha}_{x} \circ \widetilde{\alpha}_{y}(A)$. Hence $\{\widetilde{\alpha_{x}}: x \in P\}$ is a semigroup of endomorphisms. The faithfulness of $\alpha$ implies that $\widetilde{\alpha}_{0}$ is the identity map.
 
 The $\sigma$-weak continuity of $\{\widetilde{\alpha}_{x}: x \in P\}$ is an immediate consequence of Lemma \ref{boundary continuity}. This completes the proof. \hfill $\Box$
 %\newline
%Our first order of business is to show that there is no distinction between $E_{0}^{\Omega}$-semigroups and $E_{0}^{P}$-semigroups. If $\{\alpha_{x}\}_{x \in P}$ is an $E_{0}^{P}$-semigroup then $\{\alpha_{x}\}_{x \in \Omega}$ is an $E_{0}^{\Omega}$-semigroup. 
%The next lemma shows that an $E_{0}^{\Omega}$-semigroup extends to  a unique $E_{0}^{P}$-semigroup.

%\begin{lmma}
%Let $\{\alpha_{x}: x \in \Omega\}$ be an $E_{0}^{\Omega}$-semigroup on a von-Neumann algebra $M$. Assume that $\alpha$ is faithful in the sense that $\bigcap_{x \in \Omega} Ker (\alpha_{x})=\{0\}$. Then there exists a unique $E_{0}^{P}$-semigroup $\{\widetilde{\alpha}_{x}:x \in P\}$ on $M$ such that $\widetilde{\alpha}_{x}=\alpha_{x}$ for $x \in \Omega$. 
%\end{lmma}
%\textit{Proof:} 
%\newline

%\begin{rmrk}
%\%label{variant}
%Proposition \ref{technical lemma} still holds when restricted to $\Omega$. The precise statement is as follows.
%\newline
 %Let $\alpha:=\{\alpha_{x}: x \in \Omega\}$ be a semigroup of unital endomorphisms on $M$.  Assume that $\alpha$ is faithful i.e. $\displaystyle \bigcap_{x \in \Omega} Ker(\alpha_{x})=\{0\}$. 
%Suppose that $\alpha$ is $\sigma$-weakly measurable i.e. the map $\Omega \ni x \to \langle T  , \alpha_{x}(A) \rangle \in \mathbb{C}$ is measurable  for $A \in M$ and $T \in M_{*}$.
  %Then $\alpha$ is $\sigma$-weakly continuous i.e.  the map $\Omega \ni x \to \langle T, \alpha_{x}(A)\rangle \in \mathbb{C}$ is continuous for $A \in M$ and $T \in M_{*}$.
%\end{rmrk}
The following corollary is an immediate consequence of Lemma \ref{boundary continuity} and Prop. \ref{technical lemma}.
\begin{crlre}
Let $\alpha:=\{\alpha_{x}:x \in P\}$ be a $\sigma$-weakly measurable $E_{0}^{P}$-semigroup on a von-Neumann algebra $M$ with separable predual. Assume that $\alpha$ is faithful. Then $\alpha$ is an $E_0^{P}$-semigroup on $M$. 

\end{crlre}

We also need the following extension of Prop. 2.3.1 of \cite{Arveson}. 

\begin{ppsn}
\label{continuity of cocycles}
Let $\alpha:=\{\alpha_{x}: x \in P\}$ be an $E_{0}^{P}$-semigroup on $M$. Assume that $\alpha$ is faithful. Suppose $\{U_{x}:x \in \Omega\}$ is a $\sigma$-weakly measurable family of isometries in $M$ such that $U_{x}\alpha_{x}(U_{y})=U_{x+y}$ for every $x,y \in \Omega$. Then there exists a unique $\sigma$-weakly continuous family of isometries $\{\widetilde{U}_{x}:x \in P\}$  in $M$ such that  $\widetilde{U}_{x}\alpha_{x}(\widetilde{U}_{y})=\widetilde{U}_{x+y}$ and $\widetilde{U}_{x}=U_{x}$ for $x \in \Omega$. 
Moreover if $\{U_{x}:x \in \Omega\}$ is a family of unitaries then $\{\widetilde{U}_{x}:x \in P\}$ is a family of unitaries.
\end{ppsn}
\textit{Proof.} The proof is similar to that of Lemma \ref{boundary continuity} and Prop. \ref{technical lemma}.   Arveson's use of Connes' $2 \times 2$ matrix trick in Prop. 2.3.1 of \cite{Arveson} can be adapted to obtain that $\{U_{x}: x \in \Omega\}$ is $\sigma$-weakly continuous. 

Fix $x \in P$. 
Let $(x_k)$ and $(y_k)$ be sequences in $\Omega$ such that $(x_k) \to x$, $(y_k) \to x$, $U_{x_k} \to U$ and $U_{y_k} \to V$ where the convergence is in the $\sigma$-weak topology. We claim that $U=V$ and $U$ is an isometry. Consider a point $s\in \Omega$. Let $n$ tends to $\infty$ in the equation  $U_s\alpha_{s}(U_{x_k})=U_{x_k+s}$ to obtain that $U_{s}\alpha_{s}(U)=U_{x+s}$. Similarly we obtain $U_s\alpha_s(V)=U_{x+s}$. 
Calculate as follows to observe that \begin{align*}
\alpha_{s}(U^{*}U)&=\alpha_{s}(U)^{*}\alpha_{s}(U) \\
                             & = \alpha_{s}(U)^{*}U_{s}^{*}U_{s}\alpha_{s}(U) \\
                             &=U_{x+s}^{*}U_{x+s}\\
                             &=1 \\
                             &=\alpha_{s}(1).                             
\end{align*}

Since $\alpha$ is faithful, it follows that $U^{*}U=1$. Now the equation $U_{s}\alpha_{s}(U)=U_{s}\alpha_{s}(V)$ and the fact that $U_{s}$ is an isometry implies that $\alpha_{s}(U)=\alpha_{s}(V)$. But $\alpha$ is faithful. Consequently, it follows that $U=V$.
Also the equation $U_{s}\alpha_{s}(U)=U_{x+s}$ implies that if $\{U_{y}: y \in \Omega\}$ is a family of unitaries then $\alpha_{s}(UU^{*})=1=\alpha_{s}(1)$ for every $s \in \Omega$. The faithfulness of $\alpha$ implies that if $\{U_{y}: y \in \Omega\}$ is a unitary family then $UU^{*}=1$. Thus $U$ is a unitary if $U_{y}$ is a unitary for every $y \in \Omega$.
 
For $x \in P$, let $\displaystyle \widetilde{U}_{x}:=\lim_{n \to \infty}U_{x_k}$ where $(x_k)$ is any sequence in $\Omega$ converging to $x$.  By what we have shown in the proceeding paragraph and by making use of the fact that bounded sets are $\sigma$-weak compact, it follows that  $\widetilde{U}_{x}$ is well-defined.  We invite the reader to work out the details along the lines of Lemma \ref{boundary continuity} and Prop. \ref{technical lemma} to see that $\{\widetilde{U}_{x}: x \in P\}$ is the desired family of isometries/unitaries. The proof is now complete. \hfill $\Box$

%\begin{rmrk}
%The reader should compare the proofs of Lemma \ref{boundary continuity} and Prop.\ref{continuity of cocycles} with that of  Lemma 2.1 in  \cite{Sundar_symmetric}.

%\end{rmrk}

\section{Proof of  Thm. \ref{Main theorem}}
 One more technical lemma that we need is the following. This is quite standard in the theory of measurable field of 
Hilbert spaces. However we recall the proof for the reader's convenience. 

\begin{lmma}
\label{Dixmier}
Let $(X,\mathcal{B})$ be a measurable space and $\{p(x):x \in X\}$ be a $\sigma$-weakly measurable family of infinite projections in $B(\clh)$.
\begin{enumerate}
\item[(1)] Then there exists a $\sigma$-weakly measurable family of 
isometries $\{w(x):x \in X\}$ in $B(\clh)$ such that $w(x)w(x)^{*}=p(x)$. 
\item[(2)] Suppose $q$ is an infinite projection in $B(\clh)$. Then there exists a $\sigma$-weakly measurable family of partial isometries $\{w(x):x \in X\}$ in $B(\clh)$ such that 
$w(x)^{*}w(x)=q$ and $w(x)w(x)^{*}=p(x)$. 
\end{enumerate}
 \end{lmma}
\textit{Proof:} For $x \in X$, let $\clh_{x}:=\{\xi \in \clh: p(x)\xi=\xi\}$.  Then $\clh_{x}$ is an infinite dimensional Hilbert space for every $x \in X$. Let 
\[
\Gamma:=\{f:X \to \clh: ~\textrm{$f$ is weakly measurable such that $p(x)f(x)=f(x)$}\}. 
\]
 Let $\{e_1,e_2,\cdots,\}$ be an orthonormal basis for $\clh$. Observe that $\{p(.)e_{i}: i=1,2,\cdots\} \subset \Gamma$ and for every $x \in X$, $\{p(x)e_i:x \in X\}$ is total in $\clh_x$. By Prop. 7.27 of \cite{Folland}, it follows that there exists $u_1,u_2,\cdots \in \Gamma$ such that for every $x \in X$, $\{u_1(x),u_{2}(x),\cdots\}$ is an orthonormal basis for $\clh_{x}$. Let $x \in X$ be given. Let $w(x)$ be the isometry sending $e_k \to u_k(x)$. It is now clear that $\{w(x):x \in X\}$ is a $\sigma$-weakly-measurable family of isometries such that for $x \in X$, $w(x)w(x)^{*}=p(x)$. This proves $(1)$.

Suppose $q$ is an infinite projection in $B(\clh)$. Let $\{w(x):x \in X\}$ be a family of isometries as in $(1)$. Choose an isometry $v \in B(\clh)$ such that $vv^{*}=q$. It is immediate that  $\{w(x)v^{*}:x \in X\}$ is the desired family of partial isometries. This proves $(2)$ and  the proof is now complete. \hfill $\Box$.

The rest of the proof is completed as in Section 2.4 of \cite{Arveson}.  In particular, Theorem 2.4.7, Lemma 2.4.8,  Prop. 2.4.9 and Theorem 2.4.10  of \cite{Arveson} hold true for the case of a closed convex cone which is spanning and pointed. 
The proofs in \cite{Arveson} carry over to our context. However there are minor differences and a few remarks are in order.

\begin{itemize}
\item The proof of Theorem 2.4.7 in \cite{Arveson} carries over with $[0,\infty)$ replaced by $P$ and $(0,\infty)$ replaced by $\Omega$.
\item Lemma 2.4.8 of \cite{Arveson} with "strongly continuous" replaced by "strongly measurable" continues to hold for the case of a closed convex cone which is pointed and spanning.
          Instead of Lemma 10.8.7 of \cite{Dixmier}, one uses Lemma \ref{Dixmier}.
\item Proposition 2.4.9 of \cite{Arveson} for the case of a closed convex cone  is proved as in \cite{Arveson} by appealing to Prop. \ref{technical lemma}.
\end{itemize}

The proof of Thm. \ref{Main theorem} is now completed in exactly the same way as   Theorem 2.4.10 of \cite{Arveson} by  repeatedly making use of  Lemma \ref{boundary continuity}, Prop. \ref{technical lemma} and Prop. \ref{continuity of cocycles}.
 \hfill $\Box$
 
 \begin{rmrk}
 One application of Theorem \ref{Main theorem} is as follows. Let $\alpha:=\{\alpha_{x}: x \in P\}$ be an $E_{0}^{P}$-semigroup on $B(\clh)$ and let $\clk$ be a separable Hilbert space. Consider the $E_0^{P}$-semigroup $\beta:=\{\alpha_{x} \otimes 1: x \in P\}$ on $B(\clh \otimes \clk)$.  It is clear that the product systems $\mathcal{E}_{\alpha}$ and $\mathcal{E}_{\beta}$ are isomorphic. For the fibre of $\mathcal{E}_{\beta}$ at the point $x \in \Omega$ is given by
 \[
 \mathcal{E}_{\beta}(x)=\{T \otimes 1 : T \in \mathcal{E}_{\alpha}(x)\}.
 \]
 Thus, by Thm. \ref{Main theorem}, it follows that $\alpha$ and $\beta$ are cocycle conjugate.
 
 \end{rmrk}
 
\nocite{Hilgert_Neeb}
\nocite{Nica_WienerHopf}
\nocite{Renault_Muhly}
\nocite{Jean_Sundar}
\nocite{Faraut}

\bibliography{references}
 \bibliographystyle{plain}
 
 \noindent{\sc S.P. Murugan} (\texttt{spmurugan@cmi.ac.in})\\
         {\footnotesize  Chennai Mathematical Institute, \\
Siruseri, 603103, Tamilnadu.}\\[1ex]
{\sc S. Sundar}
(\texttt{sundarsobers@gmail.com})\\
         {\footnotesize  Chennai Mathematical Institute,  \\
Siruseri,  603103, Tamilnadu.}
 
\end{document}